\documentclass[submission]{arxiv2023}

\usepackage{tikz}
\usepackage{ytableau}
\usepackage{amsmath}  %% Emacs uses this line to change default reference style, I will erase this line later

\newcommand{\GL}{\operatorname{GL}}

\newcommand{\abs}[1]{\lvert #1 \rvert}
\newcommand{\Abs}[1]{\lVert #1 \rVert}
\newcommand{\bra}[1]{\langle #1 \rvert}
\newcommand{\ket}[1]{\lvert #1 \rangle}
\newcommand{\braket}[2]{\langle #1 | #2 \rangle}

\renewcommand{\Re}{\operatorname{Re}}

\newcommand{\imi}{\mathbf{i}} % imaginary i = \sqrt{-1}

\DeclareMathOperator{\Sym}{Sym}

\newcommand{\NN}{\mathbb{N}}
\newcommand{\PP}{\mathbb{P}}
\newcommand{\ZZ}{\mathbb{Z}}
\newcommand{\RR}{\mathbb{R}}

\newcommand{\CC}{\mathbb{C}}

\newcommand{\mcF}{\mathcal{F}}
\newcommand{\mcK}{\mathcal{K}}

\newcommand{\ds}{\mathrm{ds}}
\newcommand{\dt}{\mathrm{dt}}
\newcommand{\dw}{\mathrm{dw}}

\newcommand{\dz}{\mathrm{dz}}
\newcommand{\dzeta}{\mathrm{d}\zeta}

\newcommand{\domega}{\mathrm{d}\omega}

\theoremstyle{plain}
\newtheorem{thm}{Theorem}[section]

\newtheorem{prop}[thm]{Proposition}

\theoremstyle{definition}

\newtheorem{ex}[thm]{Example}

\newtheorem{rem}[thm]{Remark}

\numberwithin{equation}{section}

\definecolor{darkred}{rgb}{0.7,0,0} % darkred color
\newcommand{\defn}[1]{{\color{darkred}\emph{#1}}} % emphasis of a definition
\title[Limit shapes for skew Howe duality]{Limit shapes for skew Howe duality}

\author[D.~Betea, A.~Nazarov, and T.~Scrimshaw]{
	Dan Betea\thanks{\href{mailto:dan.betea@gmail.com}{dan.betea@gmail.com}. Supported by ERC grant COMBINEPIC No. 759702.}\addressmark{1},
	Anton Nazarov\thanks{\href{mailto:antonnaz@gmail.com}{antonnaz@gmail.com}. Partially supported by RSF grant No.~21-11-00141.}\addressmark{2},
	\and Travis Scrimshaw\thanks{\href{mailto:tcscrims@gmail.com}{tcscrims@gmail.com}. Partially supported by Grant-in-Aid for JSPS Fellows 21F51028.}\addressmark{3}}

\address{
	\addressmark{1}CNRS, LAREMA, Universit\'{e} d’Angers
	\\ \addressmark{2}Department of High Energy and Elementary Particle Physics, St.\ Petersburg State University
	\\ \addressmark{3}Department of Mathematics, Hokkaido University
}

\received{\today}

\abstract{
We study large random partitions boxed into a rectangle and coming from skew Howe duality, or alternatively from dual Schur measures.
As the sides of the rectangle go to infinity, we obtain:
1) limit shape results for the profiles generalizing the Vershik--Kerov--Logan--Shepp curve; and
2) universal edge asymptotic results for the first parts in the form of the Tracy--Widom distribution, as well as less-universal critical regime results introduced by Gravner, Tracy and Widom.
We do this for a large class of Schur parameters going beyond the Plancherel or principal specializations previously studied in the literature, parametrized by two real valued functions $f$ and~$g$.
Connections to a Bernoulli model of (last passage) percolation are explored.
}

\usepackage[backend=bibtex]{biblatex}
\usepackage[colorinlistoftodos]{todonotes}

\begin{document}

\maketitle
%% note that you DO NOT have to put your abstract here -- it is generated by \maketitle and the \abstract and \resume commands above

%%%%%%%%%%%%%%%%%%%%%%%%%%%%%%%%%%%%%%%%%%%%%%%%%%
\section{Introduction}

\paragraph{Motivation.} The study of random partitions has been, for the past 50 years, a fruitful area of research at the interface of algebraic combinatorics, representation theory, probability, and mathematical physics.
Two cornerstone results are the limit shape result of Vershik--Kerov--Logan--Shepp and the Baik--Deift--Johansson Theorem providing random matrix-like asymptotic fluctuations for the longest increasing subsequence of random permutations.
See~\cite[Ch.~1 and 2]{romik2015surprising} and references therein.

In this extended abstract we generalize both of the above by looking at a class of measures called dual Schur measures with arbitrary parameters.
Such measures fall into the category of Schur measures first introduced by Okounkov~\cite{Okounkov01}, and we use his techniques along the way.
The novelty here is twofold.
On one hand, we derive our asymptotic results keeping the two parameter sets generic and settle on a class indexed by two real functions.
On the other hand, because we are using dual Schur measures, the partitions are forced to be contained within a certain rectangle, and we see this introduces special edge behavior in the so-called critical scaling regime case (discovered by Gravner, Tracy, and Widom in~\cite{GTW01}), while allowing for universal fluctuations of the Tracy--Widom GUE type away from this regime.
As our parameters are reasonably generic, this can be seen as a form of universality.

\paragraph{Combinatorial description.} The Cauchy identity is a classical identity in algebraic combinatorics that has a repre\-sen\-ta\-tion-theoretic interpretation as characters computed in two different ways on the $\GL_n \times \GL_k$ representation $\Sym( \CC^n \boxtimes \CC^k)$.
Howe duality~\cite{howe1989remarks} is the statement that this has a multiplicity free decomposition indexed by a single partition (rather than a pair).
Using this, we can form a probability measure on partitions, a form of Okounkov's Schur measure~\cite{Okounkov01}, coming from two positive specialization of Schur polynomials $s_{\lambda}(X)$ (see~\cite[Ch.~7]{ECII} for more on them).
The Schur measure and the associated combinatorics have been well-studied from various perspectives and with a broad range of applications (\textit{e.g.},~\cite{BG15_b,Okounkov01}).
In particular, it is a quintessential example of a determinantal point process, where the joint correlation kernel is described by a determinant of pairwise correlations.

There is another $\GL_n \times \GL_k$ representation $\bigwedge (\CC^n \boxtimes \CC^k)$ that has a multiplicity free irreducible decomposition~\cite{howe1989remarks}.
This is known as \defn{skew Howe duality}.
Taking characters yields the classical dual Cauchy identity, which can be described by applying the involution $\omega$, defined by $\omega s_{\lambda}(X) = s_{\lambda'}(X)$, to the Cauchy identity.
This leads to the central object of our study, the measure on partitions given by
\begin{equation}
  \label{eq:measure-definition}
\mu_{n,k}(\lambda | X, Y) = \prod_{i,j=(1,1)}^{(n, k)} (1 + x_i y_j)^{-1} s_{\lambda}(X) s_{\lambda'}(Y),
\end{equation}
where $X = (x_1, x_2, \ldots, x_n)$ and $Y = (y_1, y_2, \ldots, y_k)$.
Now the partition $\lambda$ is constrained to a $n \times k$ rectangle.
Various specializations of this measure have appeared previously under different names: the $q$-Krawtchouk ensembles %for the principal specialization $x_i = y_i^{-1} = q^{i-1}$
(see, \textit{e.g.},~\cite{Johansson01,nazarov2022skew,nazarov2021skew}) or oriented digital boiling (see, \textit{e.g.},~\cite{GTW02,GTW01}).

\paragraph{Main results.} In this paper we study the asymptotic behavior of $\mu_{n,k}(\lambda)$ in the following sense. We first study the limit shape profile of $\lambda$ as $n, k \to \infty$ with $\lim k/n \in (0, \infty)$ and then the asymptotic behavior of the largest part $\lambda_1$ (or a closely associated quantity if $\lambda_1 = k$ by our constraints). Our main results depend on two real-valued functions $f, g \colon [0, 1] \to [0, \infty)$ --- see~\eqref{eq:param}, which we take smooth\footnote{This condition can be severely loosened, to the point we can take $f,g$ to be step functions.} and satisfying certain integrability conditions, see Section~\ref{sec:bulk} and Theorem~\ref{thm:correlation-kernel-bulk}.

\begin{thm}\label{thm:ls_intro}
  Let $k, n \to \infty$ with $k/n \to c \in (0, \infty)$ and consider parameters 
  \begin{equation} \label{eq:param}
    x_i = f(i/n), \qquad y_j = g(j/k).
  \end{equation}
  % $x_i = f(i/n), y_j = g(j/k)$
   Modulo extra technical assumptions detailed in Section~\ref{sec:bulk}, the rescaled Russian-notation profile of a random $\lambda$ sampled from $\mu_{n, k}(\cdot | X, Y)$ converges, point-wise in probability, to an explicit $1$-Lipschitz function $\Omega(t)$ supported on an explicit interval $[x_-, x_+] \subseteq [-1, c]$.\footnote{It means that $\Omega(t)$ is a straight line with slope $\pm 1$ depending on the parameters for $t$ outside $[x_-, x_+]$.}
\end{thm}

\begin{figure}
  \centering
  \includegraphics[scale=0.6]{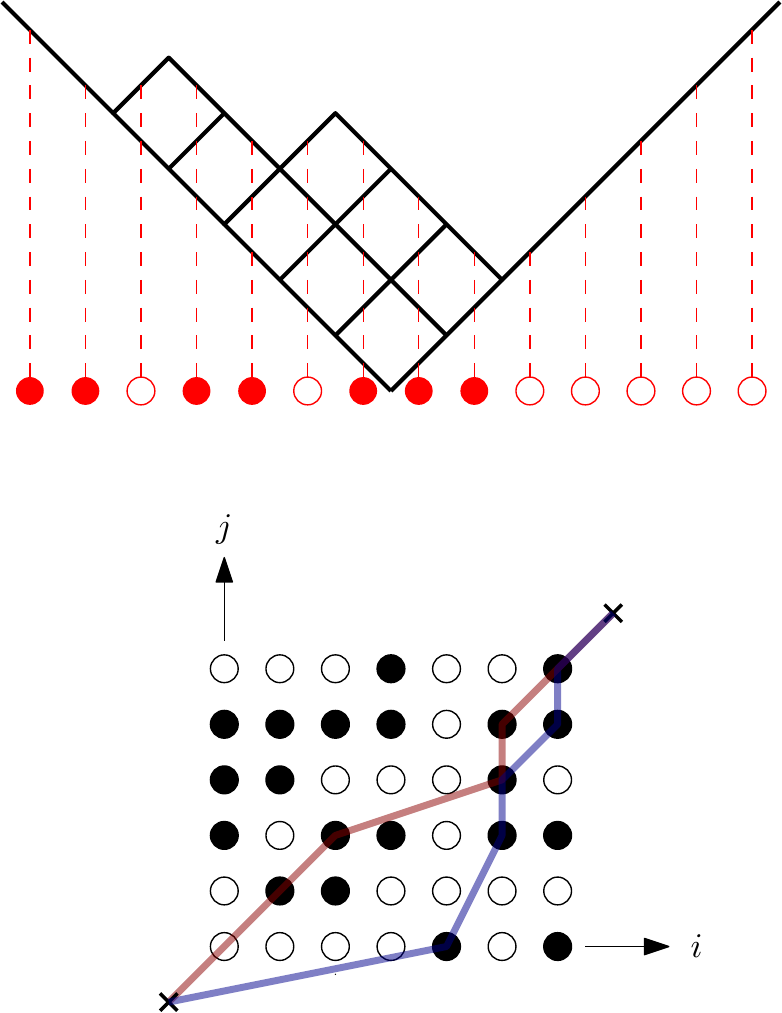} \qquad \includegraphics[scale=0.7]{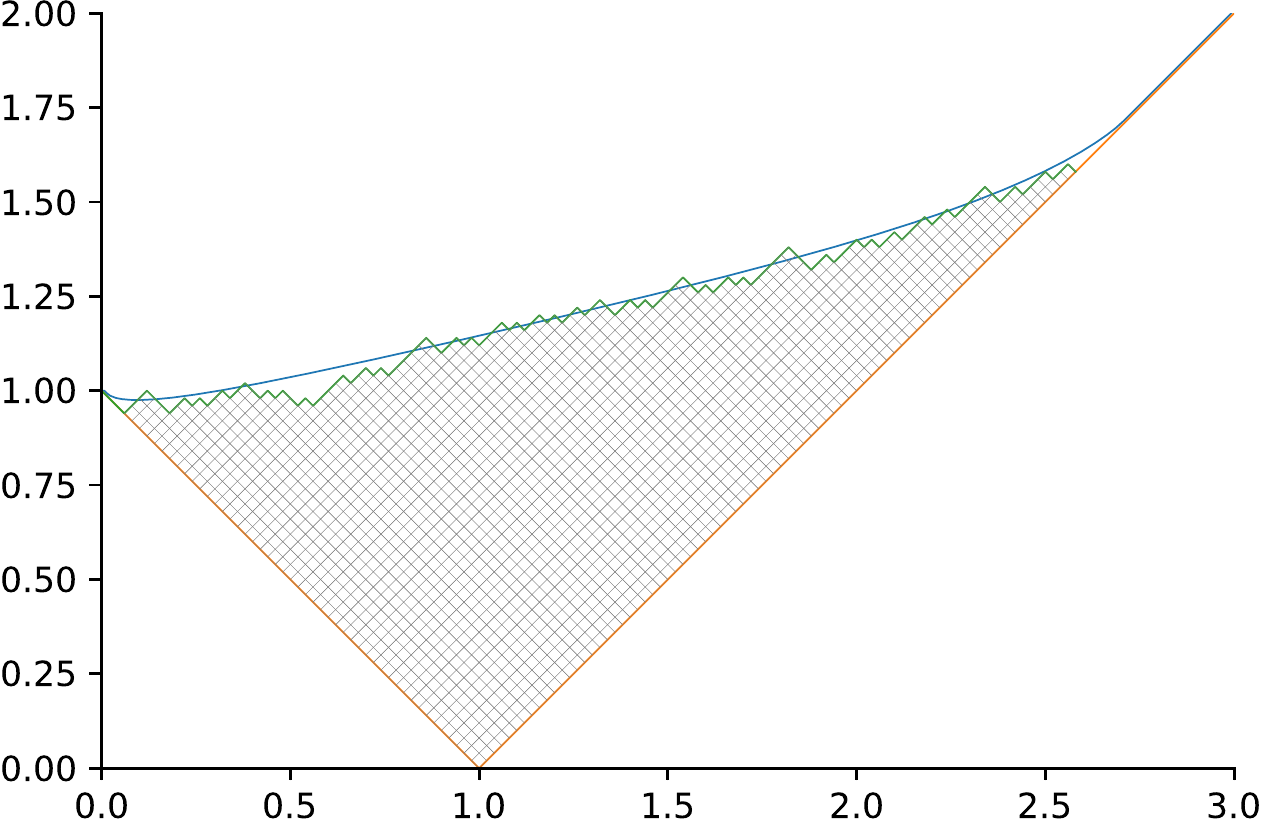}
  \caption{(Top left) Russian and Maya diagram for partition $(2^3, 1^2)$; (bottom left) two maximizing paths ($G=5$) in the Bernoulli percolation model described in the introduction (the random variables $=1$ are solid dots); (right) a random $\lambda$ from $\mu_{n,k}(\lambda|X,Y)$, the limit shape (green) superimposed ($n=50$, $k=100$, $f(x) = 3x, g(x) = x$).}
  \label{fig:maya_diagram}
\end{figure}

Our form for $X$ and $Y$ roughly means their histograms converge to continuous densities:
$n^{-1} \sum_{i} \delta_{x_i}$ converges to the continuous (with respect to Lebesgue) measure $f(s) \ds$ (and likewise for the case of the $y$'s and $g(s) \ds$) as $k, n \to \infty$ with $k/n \to c \in (0, \infty)$, where $\delta$ is the Dirac point mass.
A picture of a limit shape is depicted in Figure~\ref{fig:maya_diagram}.

We now turn to fluctuations around the limit shape.
We note that the version $g \equiv 1$ of Theorem~\ref{thm:fluct_intro} below appeared in~\cite{GTW02} and the case $f \equiv \alpha, g \equiv 1$ of Theorem~\ref{thm:crit_intro} appeared in~\cite{GTW01}, with the authors calling the latter the \defn{critical regime}.

\begin{thm}\label{thm:fluct_intro}
  With setup as in Theorem~\ref{thm:ls_intro}, consider the case when $x_+ < c$ (not full support).
  Let $L = \begin{cases} \lambda_1, &\text{if $\Omega$ convex around $x_+$}, \\ n-|\{i \mid \lambda_i = k\}|, &\text{if $\Omega$ concave around $x_+$}. \end{cases}$
  Then, with $F_{\rm GUE}$ the Tracy--Widom GUE distribution~\cite{TW94} and for some explicit constant $\sigma$, we have:
  \begin{equation}
  \lim_{n \to \infty} \PP \left( \frac{L - x_+ n}{\sigma^{-1} n^{1/3}} \right) = F_{\text{GUE}} (s).
  \end{equation}
\end{thm}

\begin{thm}\label{thm:crit_intro}
  With the scaling of Theorem~\ref{thm:ls_intro}, in the critical full support case of $x_+ = c$ and for $\Delta \in \NN$, we have
  \begin{equation}
  \lim_{n \to \infty} \PP (\lambda_1 - n c \leq -\Delta) = \det_{0 \leq i, j \leq \Delta-1} (\delta_{i, j} - K_{\rm crit}(i,j)) 
  \end{equation}
  with the matrix $K_{\rm crit} (i, j) = \sum\limits_{\ell = 0}^{(\Delta - j - 1)/2} \begin{cases} \frac{1}{2 \pi}  \frac{1}{\ell!} \sin \frac{\pi (j-i)}{2} \Gamma(\ell + \frac{j-i}{2}), & \text{if } \ell + \frac{j-i}{2} \notin \ZZ_{\leq 0}, \\ \frac{1}{2}  \frac{(-1)^\ell}{\ell! (\frac{i-j}{2} - \ell)!}, & \text{if } \ell + \frac{j-i}{2} \in \ZZ_{\leq 0}. \end{cases}$
\end{thm}

\paragraph{Connections to last passage percolation.} The statistic $\lambda_1$ for $\lambda$ from $\mu_{n, k}(\cdot|X,Y)$ has the following interpretation.
Consider the discrete $(i,j)$-grid $\{1, \dots, n\} \times \{1, \dots, k\}$. At each point $(i, j)$ place a Bernoulli $0,1$ random variable $\omega_{ij}$ distributed as $\PP(\omega_{ij} = 1) = \frac{x_i y_j}{1+x_i y_j}$.
In Figure~\ref{fig:maya_diagram} (bottom left) the non-zero $1$'s are represented as solid dots.
Let $G$ be the length of the maximizing up-right path from $(0, 0)$ to $(n+1, k+1)$\footnote{The start and end points are dummy locations for convenience.} in the following sense: the path is only allowed to go strictly up and weakly to the right, and has to maximize the number of $1$'s (solid dots in \textit{fig.\,cit.}) encountered.
Such a model is a variant of the Bernoulli last passage percolation in a random environment.
Then the Robinson--Schensted--Knuth correspondence and Schensted's theorem~\cite{ECII} yields the following.

\begin{prop}
 $G = \lambda_1$ in distribution, where $\lambda$ is distributed according to $\mu_{n, k}(\cdot|X,Y)$.
\end{prop}

%%%%%%%%%%%%%%%%%%%%%%%%%%%%%%%%%%%%%%%%%%%%%%%%%%
\section{Free fermions and finite-size correlations}
\label{sec:kernel}

We briefly recall the free fermion Fock space and half vertex operators as a method to construct the skew Schur functions as matrix elements.
Using this, we derive the correlation kernel though an application Wick's theorem.

% \begin{figure}
%   \centering
%   \includegraphics[scale=0.7]{maya_partition_2.pdf} \qquad \includegraphics[scale=0.7]{lpp_bernoulli.pdf} \quad  \includegraphics[scale=0.4]{linear-n-50-k-100-a-3-b-1.pdf}
%   \caption{Left: Russian and Maya diagram for the partition $(2^3, 1^2)$; center: two maximizing paths ($G=5$) in the Bernoulli percolation model described in the introduction (the random variables $=1$ are solid dots); right: a random partition from $\mu_{n,k}(\cdot|X,Y)$, the limit shape (green) superimposed, for $n=50$, $k=100$, $f(x) = 3x, g(x) = x$.}
%   \label{fig:maya_diagram}
% \end{figure}

% \begin{figure}
% \[
% \begin{tikzpicture}[scale=.7,>=latex]
% \draw[-,very thick] (0,0) -- ++(6.2,6.2);
% \draw[-,very thick] (0,0) -- ++(-6.2,6.2);
% \draw[<->,thick] (-6.5,0) -- (6.5,0);
% \foreach \la/\i in {4/1,3/2,2/3,2/4,2/5}
%   \draw[-] (\i,\i) -- ++(-\la,\la);
% %  \draw[-] (45:\i) -- ++(135:\la);
% \foreach \la/\i in {5/1,5/2,2/3,1/4} {
%   \draw[-] (-\i,\i) -- ++(\la,\la);
%   \draw[color=blue,dashed] (-\i,\i) + (\la+.5,\la-.5) -- (\la-\i+.5,0);
%   \fill[color=darkred] (-\i,\i) + (\la+.5,\la-.5) circle (0.1);
%   \fill[color=darkred] (\la-\i+.5,0) circle (0.1);
% %  \draw[-] (135:\i) -- ++(45:\la);
% }
% \foreach \i in {5,6} {
%   \draw[-] (-\i,\i) -- ++(.1,.1);
%   \draw[color=blue,dashed] (-\i+.5,\i-.5) -- (-\i+.5,0);
%   \fill[color=darkred] (-\i+.5,\i-.5) circle (0.1);
%   \fill[color=darkred] (-\i+.5,0) circle (0.1);
% }
% \foreach \i in {-6,-5,...,6} {
%   \draw[-] (\i,0.2) -- (\i,-0.2) node[anchor=north] {$\i$};
% }
% \end{tikzpicture}
% \]
% \caption{The construction of the Maya diagram of the partition $(5,5,2,1)$.}
% \label{fig:maya_diagram}
% \end{figure}

A \defn{partition} is an infinite sequence $\lambda = (\lambda_1, \lambda_2, \ldots)$ such that $\lambda_i \geq \lambda_{i+1} \geq 0$ for all $i \geq 1$ with only finitely many positive entries.
\defn{Fermionic Fock space} (of charge $0$) $\mcF_0$ is a span of certain vectors indexed by partitions $\lambda = (\lambda_1, \lambda_2, \dotsc)$ defined by $\ket{\lambda} = v_{a_1} \wedge v_{a_2} \wedge v_{a_3} \cdots$, where $(a_i = \lambda_i - i + \frac{1}{2} )_{i=1}^{\infty}$ is the \defn{Maya sequence} of $\lambda$, in a semi-infinite wedge space.
These vectors are obtained from the vacuum vector $\ket{0} = -\frac{1}{2} \wedge -\frac{3}{2}\wedge -\frac{5}{2}\wedge \cdots$ by the sequence of actions of the fermionic creation operators $\psi_{m}$ that adjoin $m \in\ZZ+\frac{1}{2}$ to the Maya sequence with the proper sign, and fermionic annihilation operators $\psi^*_{m}$ that remove $m$.
They satisfy the canonical anti-commutation relations $\{\psi_m, \psi^*_\ell\} = \delta_{m, \ell}$.\footnote{Here $\{a, b\} = ab+ba$, and all other anti-commutators vanish: $\{\psi_m, \psi_\ell\} = \{\psi^*_m, \psi^*_\ell\} = 0$.}
Let $\bra{\mu}$ denote the dual basis element with the natural pairing $\braket{\mu}{\lambda} = \delta_{\mu, \lambda}$.
From these operators we construct their generating series $\psi(z) = \sum_{m \in \ZZ+1/2} \psi_m z^m, \psi^*(w) = \sum_{m \in \ZZ+1/2} \psi^*_m w^{-m}$, as well as the half-vertex operators $\Gamma_{\pm}(X) = \exp \sum_{j \geq 1} \frac{p_j(x) \alpha_{\pm j}}{j}$ with $p_j(X) = \sum_i x_i^j$ and $\alpha_j = \sum_{m\in \ZZ+1/2} \psi_{m-j} \psi_m^*$.
The following properties are useful:
\begin{equation}
\label{eq:vertex_relations}
\begin{gathered}[]
[\Gamma_+(X), \Gamma_+(Y)] = [\Gamma_-(X), \Gamma_-(Y)] = 1, \quad [\Gamma_+(X), \Gamma_-(Y)] = H(X; Y) \\
	\Gamma_+(X) \ket{0} = \ket{0}, \qquad \bra{0} \Gamma_-(Y) = \bra{0}, \\
	\mathrm{Ad}_{ \Gamma_{\pm}(X) } \psi(z) = H(X; z^{\pm}) \psi(z), \quad	\mathrm{Ad}_{ \Gamma_{\pm}(X) } \psi^{\dagger}(w) = H(X; w^{\pm1})^{-1} \psi^{\dagger}(w), \\
	\bra{0} \psi(z) \psi^{\dagger}(w) \ket{0} = \sum_{\ell=0}^{\infty} z^{-1/2-\ell} w^{1/2+\ell} = \frac{\sqrt{zw}}{z - w} \quad \text{for } \abs{z} < \abs{w}
\end{gathered}
\end{equation}
where
$H(X; Y) := \prod_{i,j} \frac{1}{1 - x_i y_j}$, $[u, v] := u v u^{-1} v^{-1}$, and $\mathrm{Ad}_u x := u x u^{-1}$.
%\qquad\qquad
We denote $\Gamma'_{\pm}(Y) := \Gamma_{\pm}^{-1}(-Y)$ and
$
E(X; Y) = \prod_{i,j} (1 + x_i y_j) = H(X; -Y)^{-1}.
$

Wick's theorem allows us to write matrix elements $\bra{\mu} \Theta \ket{\lambda}$, for some special operators $\Theta$ on $\mcF_0$, as determinants. In particular, Schur polynomials (and functions) are matrix elements of the $\Gamma_\pm$ operators: $\bra{\lambda} \Gamma_-(X) \ket{0} = \bra{0} \Gamma_+(X) \ket{\lambda} \ = s_{\lambda}(X)$ and $\bra{\lambda} \Gamma'_-(X) \ket{0} = \bra{0} \Gamma'_+(Y) \ket{\lambda} = s_{\lambda'}(Y)$ (the determinants are the (dual) Jacobi--Trudi formulas for the Schur functions). Finally, using~\eqref{eq:vertex_relations}, the dual Cauchy identity becomes
\begin{align}
\sum_{\lambda} s_{\lambda}(X) s_{\lambda'}(Y) = \bra{0} \Gamma_+(X) \Gamma_-'(Y) \ket{0} = E(X; Y).
\end{align}

Furthermore, again using Wick's theorem, we obtain the following form of Okounkov's correlations for the dual Schur measure:

\begin{prop}[Correlation kernel, \cite{Okounkov01}]
  \label{prop:determinantal-ensemble-corr-kernel}
  Let $\lambda$ be random sampled from $\mu_{n,k}(\lambda | X, Y)$. Fix $a_i \in \ZZ+\frac{1}{2}$, $1 \leq i \leq p$. The associated particles (Maya sequence) form a determinantal ensemble:
  \begin{equation}
    \label{eq:measure-determinantal}
    \mu_{n,k}(\lambda \text{ has a particle at position } a_i, 1 \leq i \leq p | X, Y)=
    \det\bigl[ \mcK(a_{i},a_{j}) \bigr]_{i,j=1}^{p},  
  \end{equation}
  with the correlation kernel $\mathcal{K}(m,m')$ given by the integral representation
  \begin{equation}
    \label{eq:correlation-kernel-integral-representation}
    \mcK(m,m')=\oint\!\!\!\!\!\!\oint\limits_{\abs{w}<\abs{z}}\frac{dz}{2\pi \imi z}
    \frac{dw}{2\pi \imi w} \frac{F(z)}{F(w)} \frac{w^{m'}}{z^{m}}\frac{\sqrt{zw}}{z-w}, \quad F(z):=\prod_{i=1}^{n}\frac{1}{1-x_{i}z}\prod_{j=1}^{k}\frac{1}{1+y_{j}/z}
  \end{equation}
  % with $F(z)=\prod_{i=1}^{n}\frac{1}{1-x_{i}z}\prod_{j=1}^{k}\frac{1}{1+y_{j}/z}.$
  where both contours encircle $0$ and the $z$ contour encircles all $-y_j$'s while excluding all $1/x_i$'s.
\end{prop}

%%%%%%%%%%%%%%%%%%%%%%%%%%%%%%%%%%%%%%%%%%%%%%%%%%
\section{Asymptotics for the limit shape}
\label{sec:bulk}

We discuss how to obtain the limit shape alluded to in Theorem~\ref{thm:ls_intro}.
We look at the behavior of $\mu_{n,k}(\lambda | X, Y)$ as we take the limit $n,k\to\infty$ such that $\lim\frac{k}{n} = c \in (0, \infty)$.
The partition, drawn in Russian notation, becomes a piecewise-linear function with slope $-1$ (resp.~$1$) if a particle is in (resp.\ not in) the interval $[i, i+1]$.
Rescaling by $n^{-1}$, our partition lies in the interval $[-1, c + o(1)]$ since $k \sim nc$.
The limit density of particles turns out to be well-defined and gives the desired limit curve $\Omega$.\footnote{The convergence of the corresponding (rescaled) functions to the limit shape is a more difficult question that we do not consider here; see \textit{e.g.},~\cite[Ch. 1]{romik2015surprising} for the case of the Plancherel measure.}
Thus if we show that there exists a limiting density $\rho(t)$, then we can compute the limit shape by
\begin{equation}
\label{eq:limit_shape}
\Omega(u) = 1 + \int_{-1}^u \dt \bigl( 1 - 2 \rho(t) \bigr).
\end{equation}

We demonstrate this by showing that outside of an interval we call the support, the limit density is constant ($0$ or $1$), while inside it converges to a discrete version of the random matrix sine kernel.

We make the following assumptions on the functions $f, g$:
1) $f, g \colon [0, 1] \to [0, \infty)$ are smooth;
2) the Riemann integrals $\int_{0}^{1}\frac{\ds}{(1-f(s)z)^{2}}$, $\int_{0}^{1}\frac{\ds}{(z+g(s))^{2}}$ converge for $z \in \CC$ away from the poles; and 3) the equation
\begin{equation}
  \int_{0}^{1} \; \ds \frac{f(s)z}{(1-f(s)z)^2}
  - c \int_{0}^{1} \; \ds \frac{g(s) z}{(z+g(s))^{2}} = 0
\end{equation}
has two roots $z_{\pm}\in \mathbb{C}\cup \{\infty\}$.
Furthermore we set
\begin{equation}
  \label{eq:x-plus-x-minus}
  x_{\pm} = \int_{0}^{1} \; \ds \frac{f(s)z_{\pm}}{1 - f(s)z_{\pm}} +
  c \int_{0}^{1} \; \ds \frac{g(s)}{z_{\pm}+g(s)}. 
\end{equation}

\begin{thm}[Bulk asymptotics]
  \label{thm:correlation-kernel-bulk}
  With parameters $x_i = f(i/n), y_j = g(j/k)$, as $k, n \to \infty$ with $k/n \to c \in (0, \infty)$, and assuming the conditions 1)--3) above on $f, g$.
  Then for $t \in [x_{-},x_{+}] \subseteq [-1,c]$ and $m,m' \in \ZZ+\tfrac12$ we have
  \begin{equation}
    \label{eq:correlation-kernel-bulk-limit}
    \lim_{n\to\infty}\mathcal{K}(nt+m,nt+m') = \begin{cases}
        \dfrac{\sin\bigl( \varphi \cdot (m-m') \bigr)}{\pi(m-m')} & \text{if } m \neq m', \\
        \dfrac{\varphi}{\pi} & \text{if } m = m',
      \end{cases}
    \end{equation}
    where $\varphi(t) = \arg z(t)$ is given by the argument of the solution of
    \begin{equation}
      \label{eq:zdz-S-eq-zero-Th}
      \int_{0}^{1} \; \ds \frac{f(s)z}{1 - f(s)z} +
      c \int_{0}^{1} \; \ds \frac{g(s)}{z+g(s)} - t = 0.
    \end{equation}
    Outside of $[x_-, x_+]$, $\lim_{n\to\infty}\mathcal{K}(nt+m,nt+m') \in \{0, 1\}$.
\end{thm}

We then define $\rho(t) = \lim_{n\to\infty} \mcK(nt, nt)$ for $t \in [x_-, x_+]$; for $t < -1$ (resp.~$t > c$), we have $\rho(t) = 1$ (resp.~$\rho(t) = 0$); for $t \in [-1, x_-] \cup [x_+, c]$ the situation is more complicated as it depends on the functions and $c$, but $\rho(t)$ is still identically $0$ or $1$.

We prove Theorem~\ref{thm:correlation-kernel-bulk} by using 
$\mcK(nt,nt) = \oint\oint_{\abs{w}<\abs{z}}\frac{dz}{2\pi \imi z} \frac{dw}{2\pi \imi w} e^{n(S(z)-S(w))} \frac{\sqrt{zw}}{z-w}$
for the action $S(z)=\frac{1}{n}\ln F(z) - t \ln z$, which is asymptotically
\begin{equation}
  \label{eq:action}
  S(z) \approx-\int_{0}^{1} \; \ds
  \ln(1-f(s)z)-c\int_{0}^{1} \; \ds \ln(1+g(s)/z) - x \ln z.
\end{equation}
Only critical points of the action contribute to the integral giving $\mcK(nt, nt)$.
If we have two complex conjugate roots $z, \overline{z}$ of the equation
\begin{equation}
  \label{eq:zdz-S-eq-zero}
  z\partial_{z}S(z) = \int_{0}^{1} \ds \frac{f(s)z}{1 - f(s)z} +
  c \int_{0}^{1} \ds \frac{g(s)}{z+g(s)}-t = 0, 
\end{equation}
then
$
  \rho(t)=\frac{1}{\pi}\arg z.
$
The support of the density $[x_{-},x_{+}]$ is determined by finding the solutions $t = x_\pm$ when $z = \overline{z}$ is a double critical point of the action.
There are two such $z$'s by our assumption, called $z_{\pm}$, and obtained by solving
\begin{equation}
  \label{eq:zdz-2-S-eq-zero}
  \begin{split}
    (z\partial_{z})^{2}S(z) & = \int_{0}^{1} \; \ds \frac{f(s)z}{(1-f(s)z)^2}
    - c \int_{0}^{1} \; \ds \frac{g(s) z}{(z+g(s))^{2}} = 0.
  \end{split}
\end{equation}
Plugging into~\eqref{eq:zdz-S-eq-zero} and solving for $t$, we obtain the two ends of the support $t = x_\pm$. 

\begin{ex}[Equal parameters] \label{ex:1}
  As our first example we consider  $x_i = \alpha, y_j = 1$ (and hence $f \equiv \alpha, g \equiv 1$) for all $i,j$ with fixed $\alpha \in \RR_{>0}$.
  The measure~\eqref{eq:measure-definition} becomes 
  \begin{equation}
    \mu_{n,k}(\lambda|\alpha)=\alpha^{\abs{\lambda}} (1 + \alpha)^{-nk} \dim V_{\GL_{n}}(\lambda) \dim V_{\GL_{k}}(\lambda')
  \end{equation}
  and was considered in~\cite{GTW01} in relation to a stochastic growth process.
  We note the ensemble for $\lambda$ is Johansson's Krawtchouk orthogonal polynomial ensemble~\cite{Johansson01}.
  Substituting $f, g$ into~\eqref{eq:zdz-S-eq-zero}, \eqref{eq:zdz-2-S-eq-zero} we get $x_{\pm} = \frac{\alpha (c-1) \pm 2\sqrt{\alpha c}}{\alpha + 1}$
  and roots
  \begin{equation}
  \label{eq:z-for-const}
  z,\overline{z}=\frac{\alpha (c-1)+t(1-\alpha) \pm \sqrt{4\alpha(t+1)(t-c)+(\alpha(c-1)+t(1-\alpha))^{2}}}{2\alpha(t+1)}.
  \end{equation}
  Therefore the limit density for $t \in [x_{-},x_{+}]$ is $\rho(t) = \frac{1}{\pi}\arg z = \frac{1}{\pi}\arccos\left( \frac{\alpha (c-1) + t(1-\alpha)}{2\sqrt{\alpha(c-t)(t+1)}} \right).$
  % Not plane partitions, but certain types of square Young tableau. Not sure worth mentioning at all; there are other things that have this limit shape
  % As discussed in \cite{nazarov2021skew}, the corresponding limit shape also describes the limit shape of plane partitions~\cite{pittel2007limit}. 
\end{ex}

\begin{ex}[$q$-weights I]\label{ex:2}
  Consider now the principal specialization of characters $x_{i}=q^{i-1}, y_{j}=q^{j-1}$ related to $q$-dimensions of the corresponding representations.
  We have
  \begin{equation}
    \mu_{n,k}(\lambda|q)=\frac{q^{\Abs{\lambda}}\dim_{q}\left(V_{\GL_{n}}(\lambda)\right)\cdot q^{\Abs{\overline{\lambda}'}} \dim_{q}\left(V_{\GL_{k}}(\overline{\lambda}')\right)}{\prod_{i=1}^{n}\prod_{j=1}^{k}(q^{i-1}+q^{j-1})},
  \end{equation}
  where $\Abs{\lambda} = \sum_{i=1}^{n}(i-1)\lambda_i$ and $\overline{\lambda}$ means the complement of $\lambda$ inside the $n \times k$ rectangle.
  The $q$-dimension $\mathrm{dim}_{q}$ of the irreducible $\GL_{n}$ representation is defined in, \textit{e.g.},~\cite[\S10.10]{kac90}.
  Using explicit formulas for $q$-dimensions, one can show that the point ensemble for $\lambda$ is the so-called $q$-Krawtchouk orthogonal polynomial ensemble; see~\cite{nazarov2022skew} for details.
  To compute the limit shape, we fix $\gamma > 0$ and set $q = e^{-\gamma/n} \to 1$ as $n, k \to \infty$ (\textit{i.e.}, $f(s) = g(s) = e^{-\gamma} s$).
  The rescaled profile for a random $\lambda$ from $\mu_{n,k}(\lambda | X, Y)$ converges point-wise in probability to the limit shape given by~\eqref{eq:limit_shape} with $\rho(t)$ and $x_\pm$ given by
  \begin{align}
    \label{eq:rho-principal} 
    \rho(t) & = \dfrac{1}{\pi}
    \arccos\left(\mathrm{sgn}(-\gamma)\dfrac{e^{\gamma-\frac{\gamma  (t+1)}{2}}}{2}
      \dfrac{1-e^{\gamma(c-1)}}{\sqrt{(1-e^{\gamma (t+1)})(1-e^{\gamma(c-t)})}}\right),
    \\
      \label{eq:x-pm-principal}
     x_{\pm} & = -\frac{\mathrm{sgn}(\gamma)}{\gamma}\ln \frac{3e^{(c+1)\gamma}-e^{c\gamma}-e^{\gamma} + 3 \mp 2 \sqrt{2} \sqrt{(e^{\gamma}-1)(e^{c\gamma}-1)(e^{(c+1)\gamma}+1)}}{(1+e^{c\gamma})^{2}}.
  \end{align}
\end{ex}

\begin{ex}[$q$-weights II] \label{ex:3}
  We can also consider $x_{i} = q^{i-1}, y_j = q^{j-1}$; $\mu$ has a similar form as above:
  $\mu_{n,k}(\lambda|q,q^{-1})=\frac{q^{\Abs{\lambda}}\dim_{q}\left(V_{\GL_{n}}(\lambda)\right)\cdot q^{-\Abs{\overline{\lambda}'}}\dim_{1/q}\left(V_{\GL_{k}}(\overline{\lambda}')\right)}{\prod_{i=1}^{n}\prod_{j=1}^{k}(q^{i-1}+q^{1-j})}$.
  In the same regime as above as $q \to 1$ the limit density and support are given by
  \begin{align}
    \label{eq:rho-principal-inverse}
    \rho(t) & = \frac{1}{\pi}
    \arccos\left(\mathrm{sgn}(-\gamma)
      \frac{e^{\frac{\gamma}{2}(t+1-c)}}{2}
      \dfrac{1-e^{\gamma c}-e^{\gamma (c-t-1)}+e^{\gamma (c-t)}}
      {\sqrt{(1-e^{\gamma (t+1)})(1-e^{\gamma(c-t)})}} \right),
    \\
    \label{eq:xpm-principal-inv-principal}
    x_{\pm} & = -1-\frac{1}{\gamma}\ln 2+\frac{1}{\gamma}\ln\left(1+e^{\gamma(c+1)}\pm\sqrt{(e^{2\gamma c}-1)(e^{2\gamma}-1)}\right).
  \end{align}
\end{ex}

\begin{ex}
  Take $f(s) = \alpha s, g(s) = s$.
  We do not have closed-form formulas anymore since the equations \eqref{eq:zdz-S-eq-zero}, \eqref{eq:zdz-2-S-eq-zero} take the form: $-\frac{1}{\alpha z}\ln(1 - \alpha z) - 1 + cz \ln\left(\frac{z}{1+z}\right) + c - t = 0$, $\frac{(\alpha z - 1) \ln(1 - \alpha z) - \alpha z}{\alpha z (\alpha z - 1)} + c \left( z\ln\left(\frac{z}{z+1}\right) + \frac{z}{z + 1} \right)= 0$. Yet by solving these numerically we can still find the support $[x_-, x_+]$ and the limit shape as shown for $\alpha=3$ in Figure~\ref{fig:maya_diagram} (right).
\end{ex}

\begin{figure}[h]
  \centering
  \includegraphics[width=0.35\linewidth]{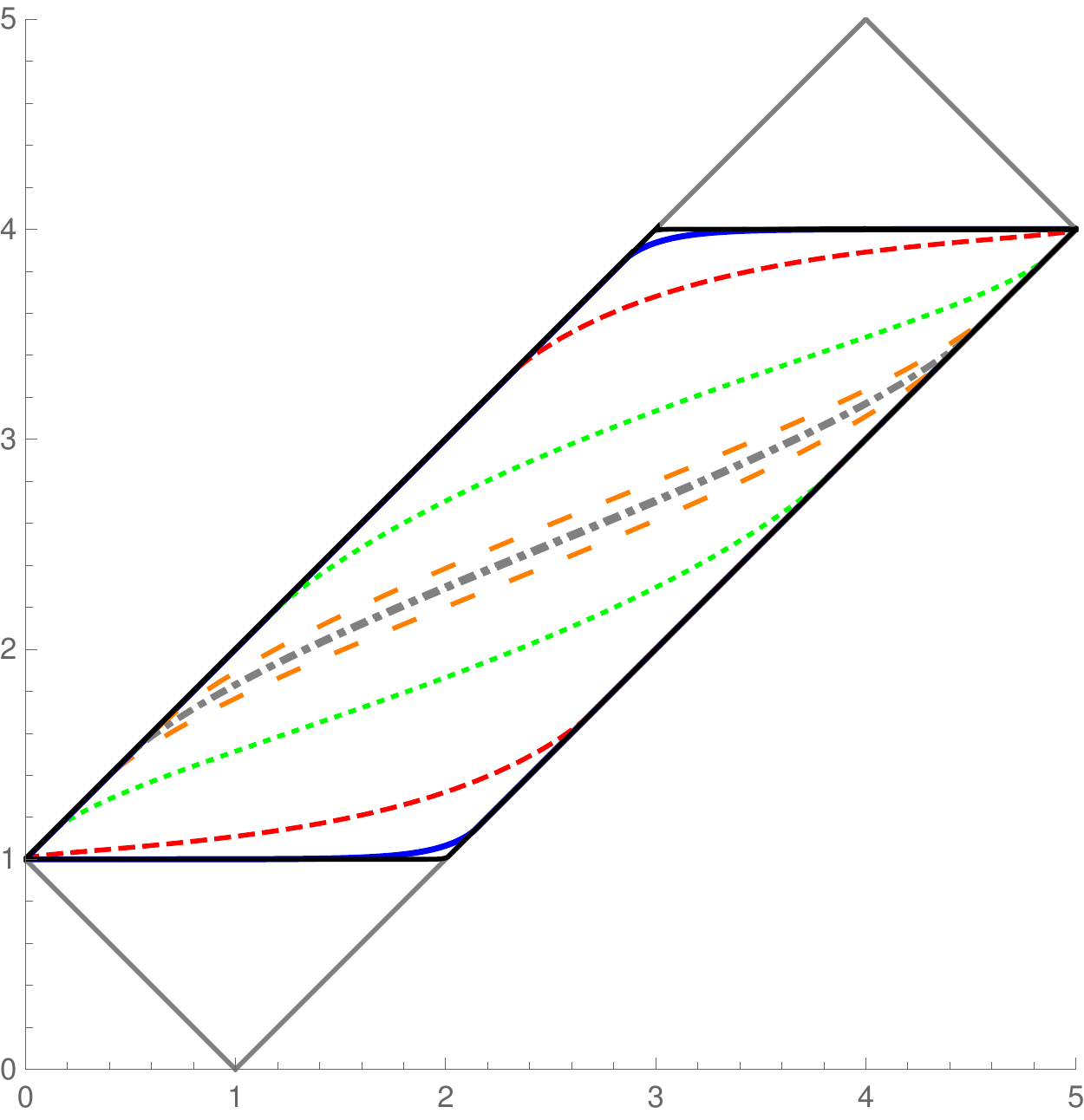}   \qquad \qquad  
  \includegraphics[width=0.35\linewidth]{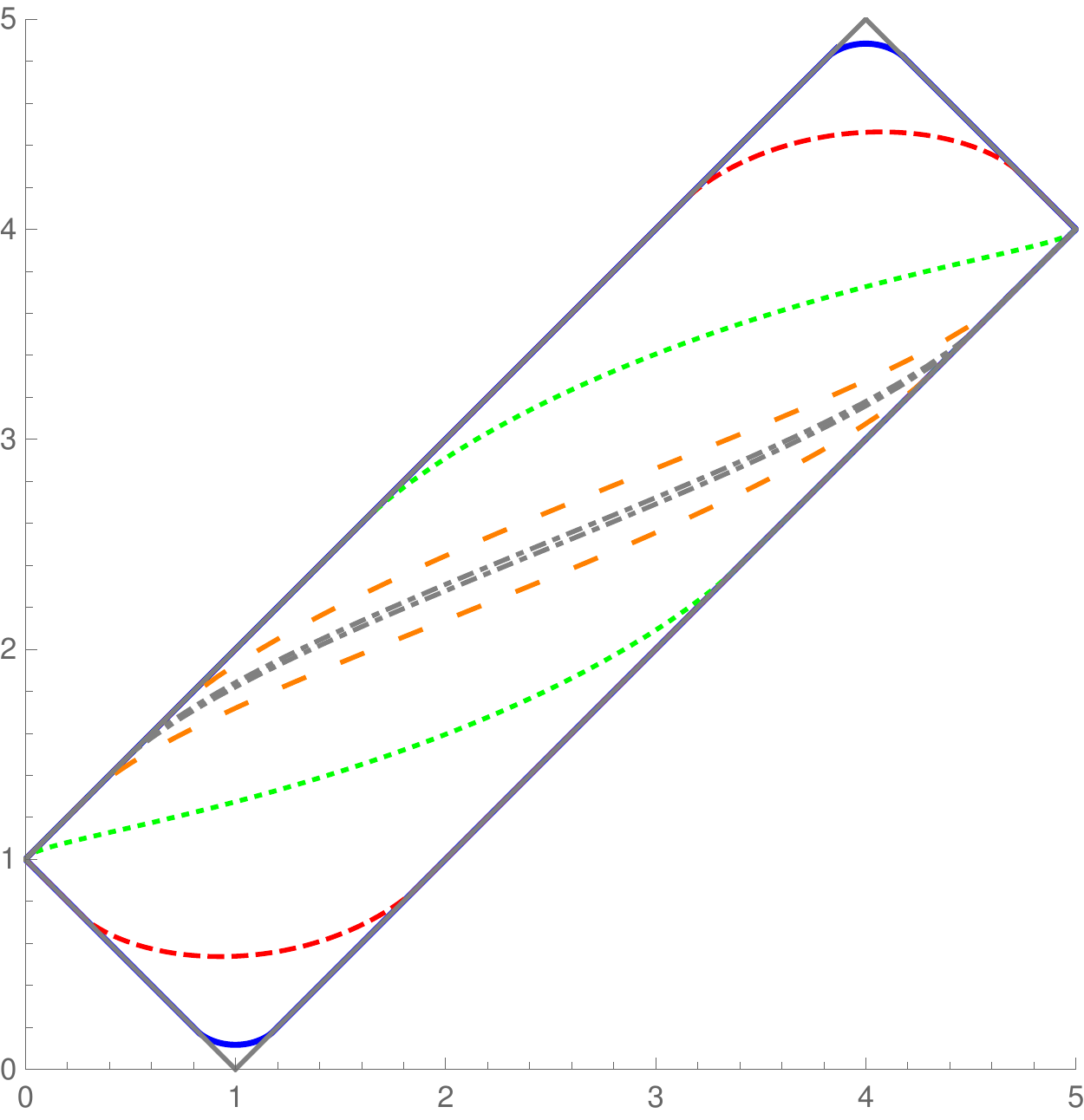}
  \caption{Plots of the limit shapes for Young diagrams corresponding
    to the densities \eqref{eq:rho-principal} (on the
    left) and \eqref{eq:rho-principal-inverse} (on the right) for $c=4$ and the
    values of $\gamma$ (bottom to top): $-10$ (solid blue), $-2$ (dashed red), $-0.5$
    (dotted green), $-0.1$ (sparsely dashed orange),  $-0.01$ (dot-dashed
    gray), $0.01$ (dot-dashed gray), $0.1$ (sparsely dashed orange), $0.5$
    (dotted green), $2$ (dashed red), $10$ (solid blue). Solid black
    lines on the left panel correspond to $\gamma=\pm\infty$
    ($q=\mathrm{const}$). We have shifted $t \mapsto t + 1$.} %% Maybe add a phrase on q<1 q>1
    % Perhaps make the $\gamma$ values into a key overlayed on the figure(s)? -- TCS
  \label{fig:limit-shape-plots}
\end{figure}

For Examples~\ref{ex:1}, \ref{ex:2} and~\ref{ex:3}, limit shapes for various values of $\gamma$ are presented in Figure~\ref{fig:limit-shape-plots}.
For $\gamma=0$ we recover the constant specialization $x_{i}=y_j=1$.

%%%%%%%%%%%%%%%%%%%%%%%%%%%%%%%%%%%%%%%%%%%%%%%%%%
\section{Boundary asymptotics}
\label{sec:boundary}

To study boundary asymptotics, we zoom in around $x_+$ and see what happens to the process locally.
In the asymptotic analysis described in the previous section, the two roots $z, \overline{z}$ coincide to a double critical point we called $z_+$. 

\paragraph{Sketch of proof of Theorem~\ref{thm:fluct_intro}.} Before we proceed further, we assume for simplicity $c$ is such that $\int_0^1 \ds \, f(s) < c \int_0^1 \ds\, \frac{1}{g(s)}$.
This means, \textit{a posteriori}, in the notation of Theorem~\ref{thm:fluct_intro}, that we are looking at fluctuations for $\lambda_1 \ll k$.
In other words, the limit shape $\Omega(t)$ is convex around $x_+$, and in turn $x_+ < c$, so we do not have full support for $\Omega$.
In Figure~\ref{fig:limit-shape-plots}, we are dealing with a limit shape \emph{sitting under the diagonal of the rectangle}. 

Near the double critical point $z_+ \in \RR - \{0\}$ we have $S(z)=S'''(z_+)(z-z_+)^{3}/6+\mathcal{O}[(z-z_+)^{4}]$.
Scaling around $z = z_+ e^{\sigma \xi n^{-1/3}} \approx z_+(1+\sigma \xi n^{-1/3}+\cdots)$ as $ n\to\infty $, with $\sigma$ a yet-to-be-determined constant, we see that we need to scale the matrix entries of $K(m, m')$ as $(m, m') \approx x_{+}n+(\xi, \eta) n^{1/3} \sigma^{-1}$. Let $\sigma^{-1} = 2^{-1/3}S'''(z_+)^{1/3} z_+ > 0$.
Then one can show
\begin{equation}
  \label{eq:airy-kernel-convergence}
  n^{1/3} \mcK (m, m') \to \mcK_{\mathrm{Airy}}(\xi,\eta) := \frac{1}{(2\pi \imi)^{2}} \int_{i \RR - \epsilon} \domega \int_{i \RR + \epsilon} \dzeta \, \frac{\exp(\zeta^3/3-\xi \zeta - \omega^3/3 + \eta \omega)}{\zeta - \omega},
\end{equation}
where $0 < \epsilon \ll 1$ and $\mcK_{\mathrm{Airy}}(\xi,\eta)$ is an integral representation of the Airy kernel~\cite{TW94}. 

\begin{figure}%[htb]
  \centering
  \includegraphics[scale=0.5]{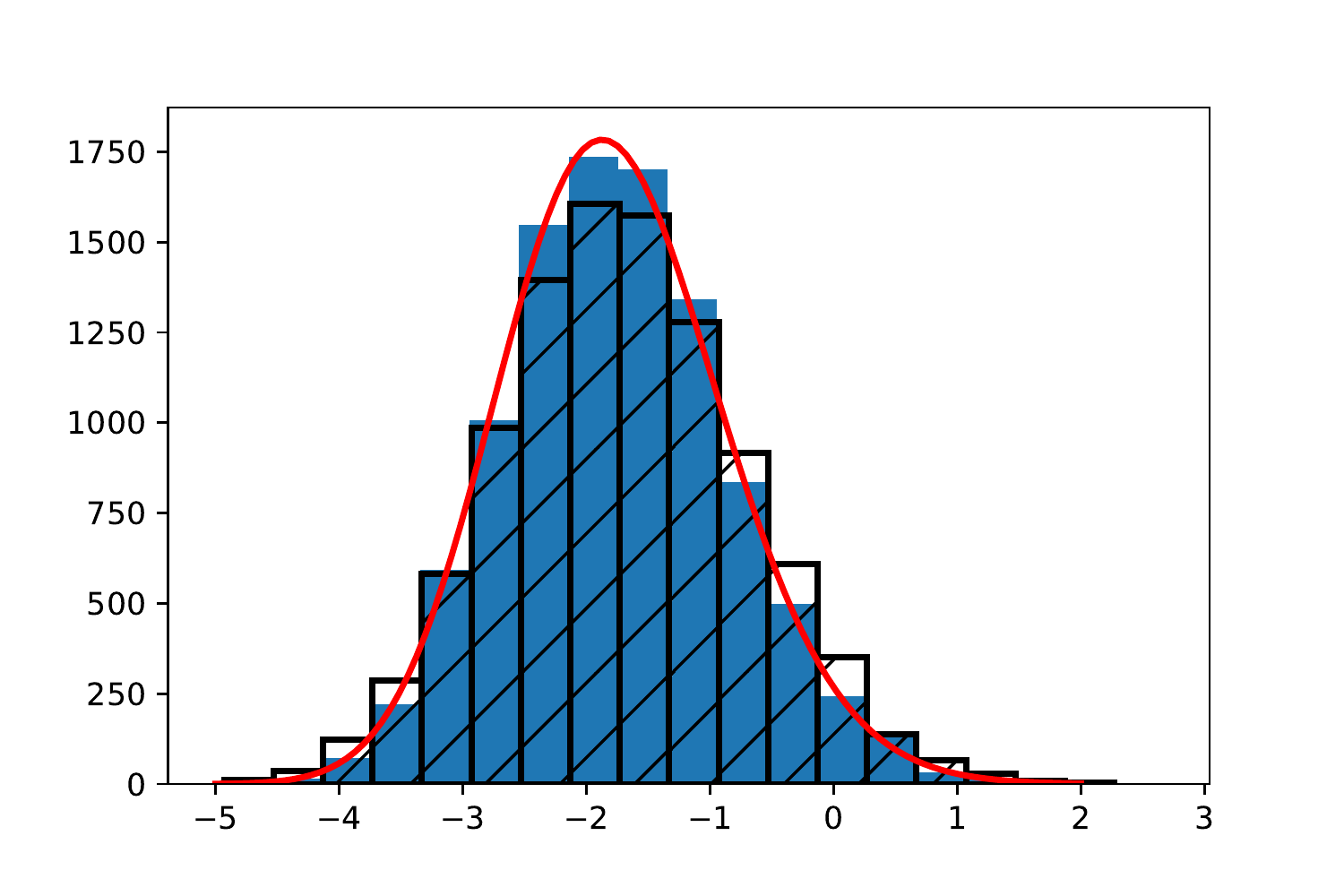}
  \includegraphics[scale=0.5]{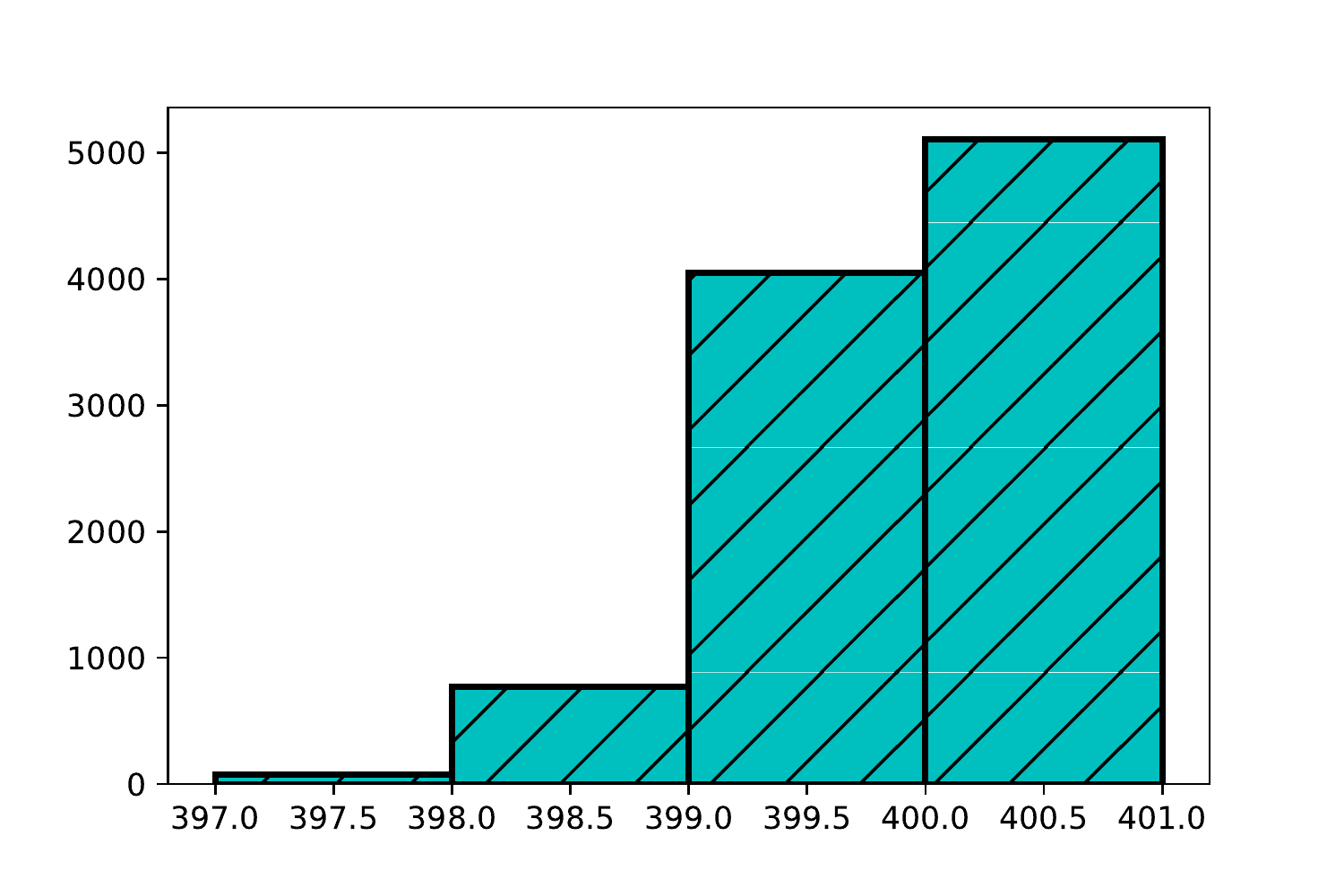}
  \caption{ (Left) Illustration of Theorem~\ref{thm:fluct_intro}: histogram of first row lengths in black ($n=200$, $k=400$, $f \equiv g \equiv 1$, $10000$ samples), Tracy--Widom distribution (red), and Tracy--Widom histogram (blue, $10000$ samples); (right) histogram of $10000$ samples of first row lengths in the critical case of Theorem~\ref{thm:crit_intro} ($n=200$, $k=400$, $f(x) = \tfrac32 x^2, g(x) = 2 x^{-1}$).
  % $x_{i}=3/2 (i/n)^{2}$, $y_{j}=2 (j/k)^{-1}$. 
  }
  \label{fig:tracy-widom-dimensions}
\end{figure}

For a quick justification notice that
\begin{equation}
n\bigl( S(z)-S(w) \bigr) \approx \frac{1}{6}\sigma^{3} z_+^{3} S'''(z_+)(\zeta^{3}-\omega^{3}) = \frac{1}{6}\sigma^{3}\left(\left.(z\partial_{z})^{3}S(z)\right|_{z=z_+}\right)(\zeta^{3}-\omega^{3})
\end{equation}
for $z=z_+ e^{\sigma \xi n^{-1/3}}, w=z_+ e^{\sigma \eta n^{-1/3}}$ in the vicinity of $z_+$, and we are therefore interested in
\begin{equation}
  \label{eq:zdz-3-S}
  (z\partial_{z})^{3}S(z)\bigg\vert_{z=z_+}=\int_{0}^{1}\ds \left(\frac{2 f^{2}(s) z_+^{2}}{(1-f(s) z_+)^{3}}+\frac{2c g(s) z_+^{2}}{(z_++g(s))^{3}}\right), 
\end{equation}
to determine $\sigma$.
Assuming this integral converges, we can conclude the convergence from~\eqref{eq:airy-kernel-convergence}.
Further technical estimates along the same lines show that Fredholm determinants converge to Fredholm determinants, and the Tracy--Widom distribution is a Fredholm determinant of the Airy kernel on $\mathcal{L}^2(s, \infty)$.
Therefore we conclude that the distribution of $\frac{\lambda_{1}-x_{+}n}{\sigma^{-1}n^{1/3}}$ is asymptotically the Tracy--Widom distribution GUE distribution. See Figure~\ref{fig:tracy-widom-dimensions} (left) for an illustration.

% If this integral diverges, the normalization constant $\sigma$ becomes infinite.
% Formally it corresponds to an infinitely-thin Tracy--Widom distribution; in practice, it means that the fluctuations are described another distribution, as we will see below.
%%% Write about integration contours here

\begin{rem}
In the case of Example~\ref{ex:1} ($f\equiv \alpha, g\equiv 1$), we have $x_{\pm}= \frac{\alpha (c-1) \pm 2\sqrt{\alpha c}}{\alpha + 1}$ and $z_+ = \frac{\alpha(c+1)-(\alpha+1)\sqrt{\alpha c}}{\alpha (\alpha c-1)}$.
This yields (recall $c, \alpha > 0$) $\sigma=\frac{(\alpha+1)c^{\frac{1}{6}}}{\alpha^{\frac{1}{6}}(\sqrt{c}-\sqrt{\alpha})^{\frac{2}{3}}(1+\sqrt{\alpha c})^{\frac{2}{3}}}$.
We see that $\sigma < \infty$ whenever $\alpha \neq c$.
\end{rem}

\begin{rem} \label{rem:crit_unif}
Continuing on the previous remark, $\alpha = c$ implies $\sigma = \infty$.
Formally this would yield an infinitely thin Tracy--Widom distribution; in practice, we obtain the critical case of Theorem~\ref{thm:crit_intro} (first discussed in~\cite[Section 3.2]{GTW01}, as the models are the same).
This happens when $\lambda_1$ hits the east-most corner of the rectangle bounding it; from the point of view of asymptotic analysis, the critical point considered is $z_+ = 0$, a singularity of the action.
\end{rem}

\paragraph{Sketch of proof of Theorem~\ref{thm:crit_intro}.} Now we turn to Theorem~\ref{thm:crit_intro}.
We land in this case when there is an obstruction to our asymptotic analysis of the generic case.
This happens when the double critical point $z_+$ of the action under consideration is $z_+ = 0$.
In this situation the constant $\sigma = \infty$, and the integral in~\eqref{eq:zdz-3-S} diverges.
Combinatorially, it happens when $\lambda_1 = nc - o(1)$ (when $\lambda_1$ just hits the east-most corner of the rectangle bounding it).
In terms of our functions $f$ and $g$, we get the following condition:
\begin{equation}
  \int_{0}^{1} \ds \, f(s) = c\int_{0}^{1}\frac{\ds}{g(s)}.
\end{equation}
Notice that $x_+=c$ is a root of~\eqref{eq:zdz-S-eq-zero} for $z=0$, and then~\eqref{eq:zdz-2-S-eq-zero} gives the above condition.

Now we cannot pass the contours through the double critical point at $0$ anymore (being a singularity of the action).
This also implies that our limiting distribution will be discrete: we will not scale the $m, m'$ coordinates continuously.
Instead we take $m = cn-h, m' = cn-h'$ for some $h,h' \in \NN+\tfrac12$ and consider the asymptotics of $K(m, m')$ as such.
The contours encircle $0$ and satisfy $\abs{w}<\abs{z}$.
The $w$ contour can be made an arbitrary small circle around $0$.
Since the contour for $z$ encircles all $-y_{j}$'s and does not contain all $1/x_{i}$'s, we can deform it to be from  $-\imi \infty$ to $\imi \infty$ with a small bump to the right of $0$ that contains contour for $w$.
We first integrate over $w$.
Since we have a critical point of the action at $0$, $S'(0)=0$, we can approximate the action as $-nS(w)\approx -nS(0) -\frac{n}{2} S''(0) w^{2}$; we thus need to compute the integral $\oint_{C_{\varepsilon}} e^{ -\frac{n}{2} S''(0) w^{2}}w^{-h'}\frac{\sqrt{zw}}{z-w}\frac{\dw}{2\pi\imi w}$.
As $\abs{w}<\abs{z}$, we can expand $\frac{1}{z-w}$ and $e^{ -\frac{n}{2} S''(0) w^{2}}$ into power series in $w$ and obtain
\begin{equation*}
  \label{eq:series-for-contour-integral}
  \oint\limits_{C_{\varepsilon}} \frac{\dw}{2\pi\imi w} e^{ -\frac{n}{2} S''(0) w^{2}}w^{-h'}\frac{\sqrt{zw}}{z-w} = \oint\limits_{C_{\varepsilon}} \frac{\dw}{2\pi\imi w} \left(\sum_{\ell=0}^{\infty} \frac{(-1)^{\ell} n^{\ell} (S''(0))^{\ell}}{2^{\ell} \ell!} w^{2\ell}\right) w^{-h'+\frac{1}{2}} \left(\sum_{j=0}^{\infty} \frac{w^{j}}{z^{j+\frac{1}{2}}}\right). 
\end{equation*}
After the residue calculation we get $\sum_{\ell=0}^{\left\lfloor h'/2-1/4\right\rfloor} \frac{(-1)^{\ell} n^{\ell} (S''(0))^{\ell}}{2^{\ell}\cdot \ell!} z^{-h'+2\ell}$.
To integrate over $z$, we note that for purely imaginary $z$ we have $\Re S(z)<0$ and therefore as $n\to\infty$ main contribution comes from the bump around $0$.
We again can approximate  $nS(z)\approx nS(0)+\frac{n}{2} S''(0) z^{2}$ and we cancel $nS(0)$ from the $w$-integral, so we need to compute the integral $\int_{-\imi \infty}^{\imi \infty} \frac{\dz}{2\pi \imi z} e^{\frac{n}{2} S''(0)z^{2}} z^{h} \left(\sum_{\ell=0}^{\left\lfloor h'/2-1/4\right\rfloor} \frac{(-1)^{\ell} n^{\ell} (S''(0))^{\ell}}{2^{\ell}\cdot \ell!} z^{-h'+2\ell} \right)$.
 We do a substitution $z=-\frac{\imi \sqrt{2t}}{\sqrt{nS''(0)}}$; the new integration contour starts at $+\infty$, loops around $0$ and then returns to $+\infty$ and we have:
 \begin{equation}
   \label{eq:corr-kernel-corner-constants}
   \mathcal{K}_{n}(h,h')= \frac{1}{2\pi}\sum_{\ell=0}^{\left\lfloor h'/2-1/4\right\rfloor}\frac{(-1)^{(h-h'+1)/2}2^{(h-h')/2-1}}{(n S''(0))^{(h-h')/2} \ell!}\int_{\infty}^{0^{+}} \dt \, t^{(h-h')/2+\ell-1} e^{-t}.
 \end{equation}
 If $(h-h')+l\geq 0$, last integral is Hankel's integral representation for the $\Gamma$-function $\Gamma(z)=-\frac{1}{2\imi \sin\pi z}\int_{\infty}^{0^{+}} \dt \, (-t)^{z-1}e^{-t}$; otherwise we obtain $ \frac{(-1)^{\frac{h-h'}{2}+l}2\pi\imi}{\left(\frac{h'-h}{2}-l\right)!}$ from the residue calculation.
 Thus the correlation kernel coincides with the formula \cite[(3.19)]{GTW01} up to a normalization.
 To find the probability $\PP (\lambda_1 - n c \leq -\Delta)$, we need to compute the determinant $\det (\delta_{i,j}-\mathcal{K}(i+1/2,j+1/2))_{i,j=0}^{\Delta-1}$ by inclusion-exclusion\footnote{We are looking at a gap probability.}; the prefactor powers with exponent proportional to $i-j$ cancel and we conclude our proof.

\acknowledgements{
The authors thank P.~Nikitin and O.~Postnova for useful discussions.
}

%% if you use biblatex then this generates the bibliography
%% if you use some other method then remove this and do it your own way
\printbibliography

\end{document}